\documentclass[11pt,a4paper,reqno]{amsart}

\usepackage{latexsym,floatflt}
\usepackage{epsfig}
\usepackage{rotating}
\usepackage{amssymb}
\usepackage[T1]{fontenc}
\usepackage{afterpage}
\usepackage{color}

\usepackage{amssymb,amsfonts,amsmath,stmaryrd}

\catcode`\@=11
\def\section{\@startsection{section}{1}%
 \z@{.7\linespacing\@plus\linespacing}{.5\linespacing}%
 {\normalfont\bfseries\scshape\centering}}
\def\subsection{\@startsection{subsection}{2}%
  \z@{.5\linespacing\@plus\linespacing}{.5\linespacing}%
  {\normalfont\bfseries\scshape}}

\def\subsubsection{\@startsection{subsubsection}{3}%
 \z@{.5\linespacing\@plus\linespacing}{-.5em}
  {\normalfont\bfseries\itshape}}
\catcode`\@=12

\newfont{\bbold}{msbm9 scaled \magstep1}
\newfont{\bbolds}{msbm7 scaled \magstep1}

\newcommand{\zs}{\mbox{\bbold Z}}

\newcommand{\bu}{\bar u}

\newcommand{\C}{\mathcal C}

\newcommand{\beq}{\begin{equation}}
\newcommand{\eeq}{\end{equation}}
\newcommand{\gf}{generating function}
\newcommand{\gfs}{generating functions}

\newcommand{\Sig}{\Sigma}

\newcommand{\Ss}{\mathcal S}

\newcommand{\G}{\mathcal G}
\newcommand{\R}{\mathcal R}
\newcommand{\ZR}{Z_\mathcal R}
\newcommand{\sR}{\tilde{\mathcal R}}
\newcommand{\cR}{{\mathcal R}^c}
\newcommand{\ZC}{Z_\mathcal C}

\newcommand{\sU}{\tilde{U}}
\newcommand{\T}{\mathcal T}

\newcommand{\cZR}{Z^c_\mathcal R}

\newcommand{\Pa}{\mathcal P}
\newcommand{\ZP}{Z_\mathcal P}

\def\bu{\bullet}
\def\ci{\circ}

\def\LL{\mathbb{L}}\def\PP{\mathbb{P}}

\def\OO{\mathbb{O}}

\def\TT{\mathbb{T}}
\def\RR{\mathbb{R}}
\def\bPP{\bar{\mathbb{P}}}

\def\ZZ{\mathbb{Z}}
\def\RR{\mathbb{R}}

\newcommand{\rrb}{\rrbracket}\newcommand{\llb}{\llbracket}

\def\la{\lambda}

\def\emm#1,{{\em #1}}

 \newtheorem{Theorem}{Theorem}
 \newtheorem{Proposition}[Theorem]{Proposition}
\newtheorem{Corollary}[Theorem]{Corollary}
\newtheorem{Lemma}[Theorem]{Lemma}

\newtheorem{Notation}{Notation}

\title[On the independence complex of square grids]
{On the independence complex of square grids}

\author{Mireille Bousquet-M\'elou}
\address{CNRS, LaBRI, Universit\'e Bordeaux 1, 351 cours de la Lib\'eration,
  33405 Talence Cedex, France}
\email{mireille.bousquet@labri.fr}

\author{Svante Linusson}
\address{Dept. of Mathematics, KTH-Royal Institute of Technology,
SE-100 44, Sweden}
\email{linusson@math.kth.se}

\author{Eran Nevo}
\address{Institute of Mathematics, Hebrew
University, Jerusalem, Israel}
\email{eranevo@math.huji.ac.il}

\date
{Mars 7, 2007}

\begin{document}

\begin{abstract}
The enumeration of independent sets of regular
 graphs is of interest in
statistical mechanics, as it corresponds to the solution of
hard-particle models. In 2004, it was  conjectured by Fendley \emm et
al., that for some rectangular grids, with toric boundary
conditions, the \emm
alternating, number of independent sets is extremely simple. More precisely,
under  a coprimality condition on the sides of the
rectangle, the number of independent sets of even
and odd cardinality always differ by $1$. 
In physics terms, this means  looking at the hard-particle model on
these grids at activity $-1$. 
This conjecture was recently proved by Jonsson.

Here we produce other families of  grid graphs, with open or cylindric
boundary
conditions, for which similar properties hold without
any size restriction: the number of independent sets of even
and odd cardinality always differ by $0, \pm 1$, or, in the cylindric
case, by some power of 2.

We show that these results reflect a stronger property  of the
independence complexes  of our graphs.  We determine the
homotopy type of these complexes  using Forman's discrete Morse
theory. We find that these complexes are either contractible,
or homotopic to a sphere, or, in the cylindric case, to a wedge of spheres.

Finally, we use our enumerative results to determine the spectra of
certain transfer matrices describing the hard-particle model on our
graphs at
activity $-1$. These results parallel certain conjectures of Fendley
\emm et al.,, proved by Jonsson in the toric case.
\end{abstract}

\maketitle

\section{Introduction}
\begin{figure}[ht]
\begin{center}
\input{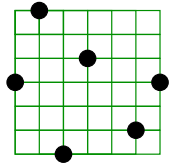}
\end{center}
\caption{A hard particle configuration --- or an independent set --- of
  the $7\times7$-grid.}
\label{fig:hardparticles}
\end{figure}

The hard-square model is a famous open problem in statistical
mechanics. In this model, some of the vertices of an $N$ by $N$ square
grid are  occupied by a particle, with the restriction that two
adjacent vertices
are never both occupied (Figure~\ref{fig:hardparticles}). In graph
theoretic terms, an admissible
configuration of particles is just an \emm independent set, of the
square grid, that is, a set of pairwise non-adjacent vertices. The key
question is to enumerate these sets by their
size, that is, to determine the following \emm partition function at
activity $u$,:
$$
Z_N(u)=\sum_{I} u^{|I|},
$$
where the sum runs over all independent sets of the grid.
This problem is highly unsolved: one does not know how to express
$Z_N(u)$, nor even the \emm thermodynamic limit, of the sequence $Z_N(u)$
(that is, the limit of $Z_N(u)^{1/N^2}$). The most natural
specialization of $Z_N(u)$, obtained for $u=1$, counts independent
sets of the $N\times N$-grid. It is also extremely
mysterious: neither the sequence $Z_N(1)$, nor the limit of
$Z_N(1)^{1/N^2}$ (the so-called \emm hard-square constant,) are known.
We refer the reader to the entry A006506 in the On-line Encyclopedia
of Integer Sequences for more details~\cite{sloane}. Note that the
thermodynamic limit of $Z_N(u)$ 
is known if one replaces the square grid by a triangular one --- a tour de
force achieved by Baxter in 1980~\cite{Baxter}.

In 2004, Fendley, Schoutens and van Eerten~\cite{Fendley} published a series of remarkable conjectures on the partition
function of the hard-square model specialized at $u=-1$. For instance,
they observed that for an $M\times N$-grid, taken with toric boundary
conditions, the partition function at $u=-1$ seemed to be equal to 1
as soon as $M$ and $N$ were coprime. They also related this conjecture
to a stronger one, dealing with the eigenvalues of the associated \emm
transfer matrices,.
These conjectures have recently been proved in a
sophisticated  way by  Jonsson~\cite{JJ}.

\begin{figure}[hb]
\begin{center}
\input{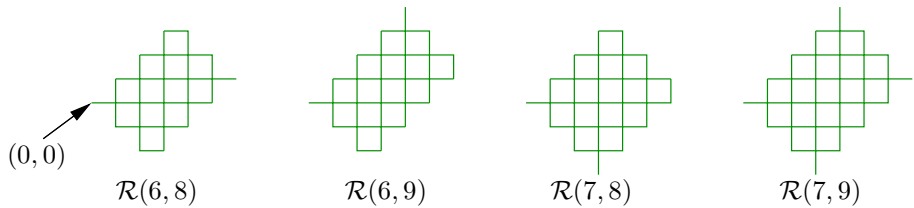}
\end{center}
\caption{The rectangular graphs $\R(M,N)$, defined in Section
\ref{sec:rectangles}.}
\label{fig:rectangles}
\end{figure}

\medskip
One of the aims of this paper is to prove that similar results hold,
 in greater generality, for other subgraphs of the square lattice,
like the (tilted) rectangles of Figure~\ref{fig:rectangles} (they
 will be defined precisely in Section~\ref{sec:rectangles}). For these
graphs, we prove that the partition function at $u=-1$ is always
$0, 1$ or $-1$. 

\begin{figure}[hbt]
\begin{center}
\input{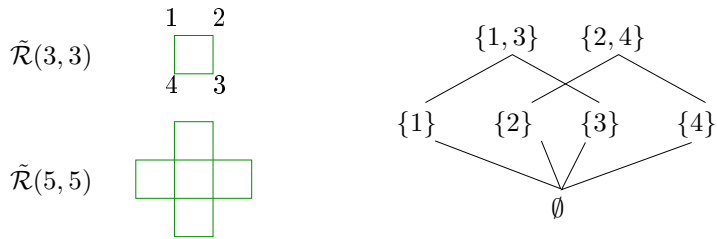}
\end{center}
\caption{Two graphs, and the independence complex of the top one. This
complex has reduced Euler characteristic 1, and is homotopic to a
0-dimensional sphere (two points). 
The patient reader can check that for the
``Swiss cross'' $\sR(5,5)$, the reduced Euler characteristic is $-1$. The
corresponding complex is homotopic to a 3-dimensional sphere.}
\label{fig:small-example}
\end{figure}

We then show that these results actually reflect a stronger
property of  the \emm independence complex, of
these graphs. The independent sets of any graph
$G$, ordered by inclusion, form a simplicial complex,
denoted by $\Sigma(G)$. See Figure~\ref{fig:small-example}, where this
complex is shown
for a $2\times2$-grid $G$. The (reduced) \emm Euler characteristic, of this
complex is, by definition:
$$
\tilde \chi_G =  \sum _{I\in \Sigma(G)} (-1)^{|I|-1}.
$$
The quantity
$|I|-1$ is the \emm dimension, of the \emm cell, $I$. The above sum is
exactly the opposite of the partition function $Z_G(u)$ 
of the hard-particle
model on $G$,  evaluated at $u=-1$. This number,
\beq\label{chi-Z}
Z_G(-1)= \sum _{I\in \Sigma(G)} (-1)^{|I|}= -\tilde \chi_G ,
\eeq
will  often be called the \emm alternating number, of independent sets.
 The simplicity of the Euler characteristic for certain graphs $G$
suggests that the complex $\Sig(G)$ could have a very simple \emm homotopy
type, (we refer to Munkres \cite{Munkres} for the topological terms
involved). We prove that this is indeed the case for various subgraphs
of the square lattice. For instance, for the rectangles of
Figure~\ref{fig:rectangles}, the independence complex is always either
contractible, or homotopy equivalent to a sphere. Our results rely on
the construction of certain \emm Morse matchings, of the complex
$\Sig(G)$.  
Roughly speaking, these matchings are parity reversing involutions on
$\Sig(G)$ having certain additional interesting properties.

Let us now  describe the contents of this paper, and compare it
to Jonsson's results.  In  Section~\ref{sec:preliminaries}, 
 we first  review   the needed background of Forman's discrete Morse theory
  \cite{Forman}.
We then describe a general
  construction of  Morse matchings for the independence complex of any
  graph. The  matchings we construct are encoded by  a \emm matching tree.,
In Sections~\ref{sec:rectangles} to~\ref{sec:paral},
we apply this general machinery to
determine the homotopy type of the independence complex of several
subgraphs of the square grid: the tilted rectangles of
Figures~\ref{fig:rectangles} and \ref{fig:smoothrectangles}
(Section~\ref{sec:rectangles}), the parallelograms of
Figure~\ref{fig:paral} (Section~\ref{sec:paral}), and variations on
them (Section~\ref{sec:quadrangles}). All these graphs
have open boundary conditions (as opposed to the toric boundary
conditions of~\cite{Fendley,JJ}). However, in Section~\ref{sec:cylindricrectangles}, we identify
two sides of the rectangles of Figure~\ref{fig:rectangles} to obtain
rectangles with cylindric boundary conditions. Again, we determine the
homotopy type of the associated independence complex. Note that Jonsson
recently went one step further by studying the same tilted rectangles
with \emm toric, boundary conditions; but he was only able to
determine the Euler characteristic
(\cite{JJdiag}, Section 7).  Our results
deal with a finer invariant of the complex (the homotopy type) and the
proofs are simpler, \emph{but} they also
refer to easier graphs; the toric case is at the moment
 beyond reach of our methods.

Finally, in Section \ref{sec:transfer} we give background about
transfer matrices, and show how the results of the previous sections can be
used to derive the spectrum, or at least part of the spectrum, of
several transfer matrices naturally associated to our graphs.

We conclude the paper with a discussion on possible extensions of our
work, with the double objective of discovering new subgraphs with a
simple alternating number of independence set, and addressing the
conjectures of Fendley \emm et al., in the cylindric case.

\section{Morse matchings for independence  complexes}
\label{sec:preliminaries}

The key tool in the proof of our results is the construction of
\emm Morse matchings, on independence complexes. In this section we 
give the background and explain how we use a matching tree as a systematic
construction for Morse matchings.
\subsection{Generalities}
Let us first recall that the poset of independent sets of a graph $G$,
ordered by inclusion, is
a simplicial complex, denoted $\Sigma(G)$.

We regard any  finite poset
$P$ as a directed graph, by considering the Hasse
diagram of $P$ with edges
pointing down (that is, from large elements to small elements).

A set $M$
of pairwise disjoint edges of this graph is called a
\emph{matching} of $P$. This matching is \emm perfect, if it covers
all elements of $P$.
This matching is \emph{Morse} (or \emph{acyclic}) if the directed
graph obtained from $P$ by reversing the direction of the edges in $M$
is acyclic. For instance, the matching of the complex of
Figure~\ref{fig:small-example} formed of the edges $(\emptyset,\{2\})$,
$(\{3\},\{13\})$, $(\{4\},\{24\})$ is acyclic.

\begin{Theorem}[\cite{FormanUser}, Theorem 6.3]\label{thm:Forman}
Let $\Sig$ be a finite simplicial complex, seen as a poset,
and $M$ a Morse matching on $\Sig$,
such that the element $\emptyset$ of $\Sig$ is matched. For $i\ge 0$, let $n_i$ be
the number of unmatched $i$-dimensional elements of $\Sig$.
Then there exists a CW-complex having $1+n_0$ $0$-dimensional cells and $n_i$
$i$-dimensional cells for $i\ge 1$ that is homotopy equivalent to $\Sig$.
\end{Theorem}
Again, we refer to~\cite{Munkres} for the topological terms involved.
We will only use the following  immediate corollary.

\begin{Corollary}\label{cor:Forman}
Under the assumptions of Theorem~{\rm\ref{thm:Forman}}, if all
unmatched elements in $\Sig$ have the same dimension $i>0$  and there
are $j$ of them,
 $\Sig$ is homotopy equivalent to a wedge of $j$  spheres of dimension $i$.
In particular, if $\Sig$ is perfectly matched, then it is contractible.
\end{Corollary}

\subsection{Matching trees}
Let us now describe the general principle  we use to  construct Morse
matchings for a complex $\Sig(G)$. Let $V$ denote the vertex set of $G$.
The most naive way to define a matching of $\Sigma\equiv \Sig(G)$ is
probably the following. Take a vertex $p \in V$, and denote by $N(p)$
the set of its neighbours. Define
$$
\Delta= \{I\in \Sig: I\cap N(p)=\emptyset\}.
$$
The set of pairs $(I,
I\cup\{p\})$, for $I\in \Delta$ and $p\not \in I$, forms a perfect
matching of $\Delta$, and hence a matching of $\Sig$. We call $p$ the
\emm pivot, of this matching. The unmatched
elements of $\Sig$ are those containing at least one element of
$N(p)$. There may be many unmatched elements, but we can now choose
another pivot $p'$  to match some elements
of $\Sig\setminus \Delta$, and repeat this operation as long as we
can. Of course, the resulting matching will depend on the successive
choices of pivots.

This rather naive idea is the leading thread in the construction of
our Morse matchings of $\Sig$. In some occasions, we will have to split
the set of yet unmatched elements into two subsets, and choose
a different pivot for each of them. This explains why our matching
procedure will be encoded by a  branching structure, namely a
\emm plane rooted tree,, called a \emm matching tree, of $\Sig$.
The nodes of this tree represent sets of yet unmatched elements. Some
nodes are
 reduced to the empty set, and all the others 
are subsets of $\Sig$ of the form
$$
\Sig({A,B}) = \{I \in \Sig : A\subseteq I \hbox { and } B \cap I =
\emptyset\},
$$
where
\beq \label{AB-conditions}
A\cap B=\emptyset \quad \hbox{ and } \quad N(A):= \cup_{a\in A}N(a) \subseteq B.
\eeq
We say that the vertices of
$A\cup B$ are the \emm prescribed vertices, of $\Sig(A,B)$. The
root of the tree is $\Sig(\emptyset, \emptyset)= \Sig$, the set of all
independent sets of $G$.
The sets $A$ and $B$ will increase along branches, making the sets
$\Sig(A,B)$ of unmatched elements smaller and smaller.  The leaves
of the tree will have cardinality 0 or 1, and will contain the
elements that are left unmatched at the end of the procedure.

Consider a node of the tree of the form $\Sig(A,B)$. How can we match
its elements? If the node has cardinality 1, that is, if $A\cup B = V$,
then we are stuck, as there is no non-trivial matching of a graph
reduced to one vertex.  If  $A\cup B
\subsetneq V$, we may match some of the elements of $\Sig(A,B)$.
Pick a vertex  $p$ in
$V':=V\setminus (A\cup B)$. Because of~\eqref{AB-conditions}, the
neighbours of $p$ are either in $B$, or in $V'$. This makes $p$ a good
\emm tentative pivot,. If we actually use $p$ as a pivot to match
elements of $\Sig(A,B)$, we will be left with the following set of
unmatched elements:
$$
U=\{I \in \Sig : A\subseteq I, I\cap B=\emptyset,  I\cap N(p) \not =
\emptyset\}.
$$
If $p$ has no neighbour in $V'$, the above set is empty, and we have
perfectly matched $\Sig(A,B)$. If $p$ has exactly one neighbour in
$V'$, say $v$, then $U=\Sig({A\cup\{v\},
  B\cup N(v)})$.
 However, if $p$ has at least two neighbours in $V'$, say $v$ and
$v'$, the set $U$ is  \emm not, of the form $\Sig(A',B')$. Indeed,
 some of the unmatched sets $I$ contain $v$, some others
don't, but then they have to contain $v'$. This puts us into trouble,
as we want to handle only unmatched sets of the form $\Sig(A',B')$. 
We  circumvent this difficulty by splitting the original set $\Sig(A,B)$
into two disjoint subsets of the form $\Sig(A',B')$, that differ by
the status of, say, the vertex $v$. More precisely, we write:
$$
\Sig(A,B)= \Sig({A,B\cup\{v\}}) \ \uplus\  \Sig({A\cup\{v\},  B\cup N(v)}),
$$
and then study separately each subset.

The above discussion justifies the following construction of the
children of a node. If
this node is the empty set (no unmatched elements), we declare it a
leaf.  Otherwise, the node is of the form
$\Sig(A,B)$. If $A \cup B=V$, the node has cardinality 1, and we also
declare it a leaf. We are left with nodes of the form
$\Sig(A,B)$, with $A\cup B
\subsetneq V$. 
Choose a vertex  $p$ (the \emm tentative pivot,\/) in
$V'=V\setminus (A\cup B)$, and   proceed as follows:
\begin{itemize}
\item
If $p$ has at most  one neighbour in $V'$,
define $\Delta(A,B,p)$ to be the subset of
 $\Sig(A,B)$ formed of sets that do not intersect $N(p)$:
$$
\Delta(A,B,p)=  \{I \in \Sig : A\subseteq I \hbox { and } B \cap I =
I\cap N(p)= \emptyset\}.
$$
Let $M(A,B,p)$ be the perfect matching of $\Delta(A,B,p)$ obtained by
 using $p$ as a
pivot. Give to the node  $\Sig(A,B)$ a unique child, namely the set
 $U=\Sig(A,B)\setminus \Delta(A,B,p)$ of unmatched elements. This set
 is empty if $p$ has no neighbour in $V'$. In this case, we say that
 $p$ is a \emm free vertex, of $\Sig(A,B)$. If $p$ has exactly one
 neighbour $v$ in $V'$, then $U=\Sig({A\cup\{v\},
  B\cup N(v)})$.  Index the new edge by the pivot $p$. 
We say that the 3-tuple $(A,B,p)$ is a \emm  matching site, of the tree.

\item Otherwise,  let us choose one neighbour $v$ of $p$ in $V'$. The node
  $\Sig(A,B)$ has two children, only differing by the status of
  $v$. More precisely, the left child is $\Sig({A,B\cup\{v\}})$
and the right child is $\Sig({A\cup\{v\},  B\cup N(v)})$.
The union of these two sets is  $\Sig(A,B)$.
Label the two new edges by the \emm  splitting vertex, $v$. 
We say that $(A,B,v)$ is a \emm splitting site, of the tree.
\end{itemize}
Observe that the new nodes satisfy Conditions~\eqref{AB-conditions},
unless they are empty.

Given a sequence of choices of tentative pivots and
splitting vertices, we obtain a matching $M$ of $\Sig$ by taking the
union of all partial matchings $M(A,B,p)$
performed at the matching sites of the tree. The unmatched elements
are those sitting  at the leaves of the tree.

The above construction is rather natural, and we invite the
reader to practice with the example given in
Figure~\ref{fig:tree-ex}. In this figure, every (non-empty) node is
described by the
vertices of $A$ (in black) and $B$ (in white). At the matching site
$(\emptyset, \emptyset,1)$, the elements $\emptyset$ and $\{1\}$
(among others) are matched. At the matching site
$(\{2\},\{1,3,4,6\},5)$, the elements $\{2\}$ and $\{2,5\}$ get
matched, among others.
At the end of the matching procedure, 
the  independent set $\{2,7\}$
is the only unmatched element of $\Sig(G)$.

\begin{figure}[htb]
\begin{center}
\input{tree-example.pstex_t}
\end{center}
\caption{A subgraph $G$ of the square grid, and one of its matching trees.}
\label{fig:tree-ex}
\end{figure}

We now aim at showing
that the matchings obtained with the above procedure are in fact Morse.
The following lemma gathers some properties of this construction.
 \begin{Lemma}\label{lemma:f-is-Morse} Every matching tree
   satisfies the following properties:
   \begin{enumerate}
     \item For every matching site $(A,B,p)$, the matching $M(A,B,p)$
     is a Morse matching of $\Delta(A,B,p)$ (still ordered by inclusion).
     \item
      Let $(A,B,p)$ be a matching site with a non-empty child
      $\Sig({A\cup\{v\}, B\cup N(v)})$. Let $I\in\Delta(A,B,p)$ and  
      $J\in\Sig({A\cup\{v\}, B\cup N(v)})$. Then $J\nsubseteq I$.
     \item
      Let $(A,B,v)$  be a splitting site, $I\in\Sig({A,B\cup\{v\}})$ and $J\in\Sig({A\cup\{v\},  B\cup N(v)})$.  Then $J\nsubseteq I$.
   \end{enumerate}

 \end{Lemma}
 \begin{proof}
   (1) Consider the Hasse diagram of the poset $\Delta(A,B,p)$ and
       its directed version, with all edges pointing down. Now, reverse
       the edges  of $M(A,B,p)$. The up edges join two elements of the
       form $I\setminus \{p\}$, $I$, so they correspond to adding
       the vertex $p$.  The down edges correspond to deleting a vertex
       different from $p$. Clearly there cannot be a directed cycle in
       $\Delta(A,B,p)$.

(2) The set $J$ contains $v$, a neighbour of the pivot $p$, while none
    of the matched sets $I$ of $\Delta(A,B,p)$ contains $v$.

(3) Here again, $J$ contains $v$, but $I$ doesn't.
 \end{proof}

The following easy lemma appears as Lemma 4.3 in Jonsson's thesis
\cite{JJthesis}.
 \begin{Lemma}\label{lemma:JJ-Morse}
Let $V$ be a finite set 
and $\Delta=\Delta_1\uplus\Delta_2$ a collection of subsets of $V$,
ordered by inclusion. 
Assume  that if $\sigma\in
\Delta_1$ and $\tau\in \Delta_2$ then $\tau\nsubseteq\sigma$. Then
the union of two acyclic matchings on $\Delta_1$ and $\Delta_2$
respectively is an acyclic matching on $\Delta$.
 \end{Lemma}
We can now establish the main result of this section.
 \begin{Proposition}\label{cor:f-is-Morse}
For any graph $G$ and any matching tree of $G$,
the matching of $\Sig(G)$ obtained by taking the
union of all partial matchings $M(A,B,p)$ performed at the matching
sites is Morse.
 \end{Proposition}
 \begin{proof}
We will prove  by backwards induction, from the leaves to the root, the
following property:
\begin{quote}
For every node $\tau$ of the matching tree, the union of the partial matchings
performed at the descendants of $\tau $ (including $\tau $ itself) is
a Morse matching of $\tau$. We denote this matching $U\!M(\tau)$ (for
 Union of Matchings).
\end{quote}
The leaves of the tree are either empty sets or singletons, endowed
with the empty matching, which is Morse. This supplies the induction
base. Consider now a non-leaf node of the tree, of the form
$\tau=\Sig(A,B)$.

Assume $(A,B,p)$ is a matching site. By Lemma
\ref{lemma:f-is-Morse}.1, $M(A,B,p)$ is Morse.
 If the (unique) child of $\tau$ is empty, then
 $U\!M(\tau)=M(A,B,p)$ and we are done. If this
child is  $\tau'=\Sig({A\cup\{v\},B\cup N(v)})$,  the induction
hypothesis tells us that $U\!M(\tau')$ is Morse.  By Lemma
\ref{lemma:f-is-Morse}.2, we can apply
 Lemma~\ref{lemma:JJ-Morse} with $\Delta_1=\Delta(A,B,p)$,
$\Delta_2= \tau'$ and $\Delta=\tau$, where the partial matchings on
$\Delta_1$ and $\Delta_2$ are respectively $M(A,B,p)$ and
$U\!M(\tau')$. This shows that $U\!M(\tau)$ is Morse.

Assume $(A,B,v)$ is a splitting site.  By induction hypothesis we already
have Morse matchings on both children of $\tau$, namely
$\Delta_1=\Sig({A,B\cup\{v\}})$
and  $\Delta_2=\Sig({A\cup\{v\},  B\cup N(v)})$.  Again,   Lemma
\ref{lemma:f-is-Morse}.3 allows us to apply Lemma \ref{lemma:JJ-Morse}, and
this shows that the union $U\!M(\tau)$ of $U\!M(\Delta_1)$ and
$U\!M(\Delta_2)$ 
is  Morse.

This completes the induction. The case where $\tau$ is the root of the
tree gives the proposition.
 \end{proof}

\section{Rectangles with open boundary conditions}
\label{sec:rectangles}

In what follows, we consider $\zs^2$ as an infinite graph, with edges
joining vertices at distance~1 from each other. For $M, N \ge 1$, let
$\R(M,N)$ be the subgraph of $\zs^2$ induced by the points $(x,y)$ satisfying
$$
y\le x\le y+M-1 \quad  \mbox{ and } \quad -y\le x\le -y+N-1.
$$
Examples are shown on Figure~\ref{fig:rectangles}.
Note that $\R(M,N)$ contains $\lceil \frac{MN}{2}\rceil$ vertices.
Other rectangular shapes arise when we look
at the subgraph $\sR(M,N)$ of $\zs^2$ induced by the points $(x,y)$
such that
$$
y\le x\le y+M-1 \quad  \mbox{ and } \quad -y+1\le x\le -y+N
$$
(see Figure~\ref{fig:smoothrectangles}). More precisely, the graphs
$\sR(2M+1,2N+1)$  are not isomorphic to any of the $\R(K,L)$ (the
other graphs $\sR(M,N)$ \emm are, isomorphic to an $\R(K,L)$).

\begin{figure}[hbt]
\begin{center}
\input{smoothrectangles.pstex_t}
\end{center}
\caption{The graphs $\sR(M,N)$.}
\label{fig:smoothrectangles}
\end{figure}

We study the independence complexes of the rectangles $\R(M,N)$ and
$\sR(M,N)$.
Recall the general connection~\eqref{chi-Z} between the  (reduced) Euler
characteristic  of these complexes and the alternating
number of independent sets.
We use below the notation $\ZR(M,N)$ rather than $Z_{\R(M,N)}$.

\begin{Theorem}\label{thm:open}
Let $M, N \ge 1$. Denote
 $m= \lceil M/3\rceil$ and  $n= \lceil  N/3\rceil$.
\begin{itemize}
  \item
If $M\equiv_3 1 $ or $N\equiv_3 1 $, then $\Sigma(\R(M,N))$ is
  contractible and $\ZR(M,N)=0$.

 \medskip
\item
Otherwise,
   $\Sigma(\R(M,N))$ is homotopy equivalent
  to a sphere of dimension $mn-1$, and $\ZR(M,N)=(-1)^{mn}$.
\end{itemize}
The above holds also when replacing $\R(M,N)$ by $\sR(M,N)$.
\end{Theorem}
\noindent
{\bf Remark.} We will show in Corollary~\ref{coro-0} that, for
 $N\equiv_3 1$ and $M> 2^{1+\lceil N/2\rceil}$, the alternating number
 $\ZR(M,N;C,D)$ 
of independent sets of $\R(M,N)$ having border conditions $C$ and $D$
 on the two 
extreme diagonals of slope 1 is actually 0 for all configurations $C$
 and $D$.  This will indirectly follow from the
study of tilted rectangles with cylindric boundary conditions
performed in Section~\ref{sec:cylindricrectangles}.
%
\begin{proof}
We study the graphs $\R(M,N)$ and $\sR(M,N)$ together.
We  construct  Morse matchings of the independence complexes of
these graphs by
following the general principles of Section~\ref{sec:preliminaries}.
 We  need to specify our choice of tentative pivots and
splitting vertices. Our objective is to  minimize the combinatorial
explosion of cases, that
is, the number of splitting sites. Consider a node $\Sig(A,B)$.
\begin{enumerate}
\item  In general, the tentative pivot
 $p=(i_0,j_0)$ is chosen in $V'=V\setminus (A\cup B)$ so  as to minimize
 the pair $(i+j,i)$ for the lexicographic order. That  is, $p$  lies
as high as possible on the leftmost diagonal of slope   $-1$. However,
 if there is at least one free vertex (that is, a
 vertex   of $V'=V\setminus (A\cup B)$ having no neighbour in
 $V'$) on the \emm next, diagonal  $i+j=i_0+j_0+1$, then we choose one
 of them as the   pivot: the only child of  $\Sig(A,B)$ is then the empty set. 
\item
If the tentative pivot has several  neighbours in $V'=V\setminus (A\cup
B)$,  then it has exactly two neighbours in $V'$, namely its North and East
neighbours. Indeed, the other two neighbours are smaller than $p$ for the
lexicographic order, and thus belong to $A \cup B$.
Take $v$, the splitting vertex, to be the East neighbour of $p$.
\end{enumerate}
These conventions are illustrated by an example in
Figure~\ref{fig:match-open}. As before, the elements of $A$ and $B$
are indicated by $\bullet$ and $\circ$  respectively. 
Rather than labeling the vertices of the graph 
and the edges of the tree, 
we have indicated the tentative 
pivots by $*$, and the splitting vertices by $\triangle$.

\begin{figure}[thb]
\begin{center}
 \scalebox{0.9}{\input{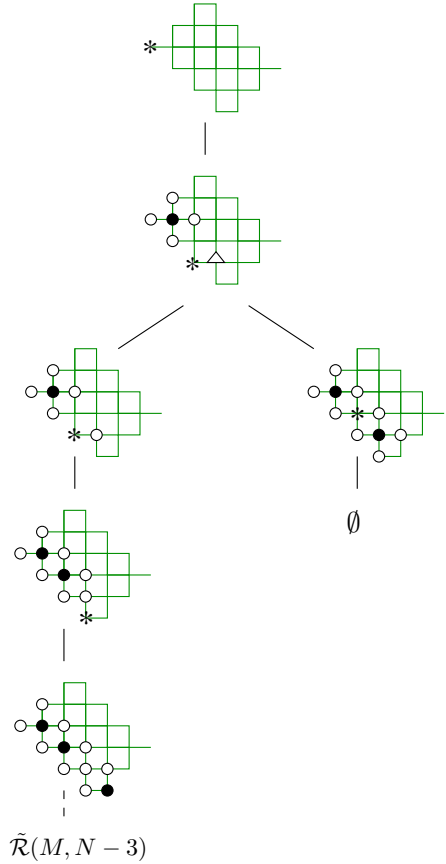}}
\end{center}
\caption{The top of the tree describing the matching of
  $\Sigma(\R(8,6))$.}
\label{fig:match-open}
\end{figure}

We are going to prove by induction on $N$ and $M$ the following properties,
valid both for the graphs $\R(M,N)$ and $\sR(M,N)$:
\begin{itemize}
  \item[$(A)$] if $M\equiv_3 1$ or $N\equiv_3 1$, then there is no unmatched element
in $\Sig$,
\item[$(B)$] otherwise, there is a unique unmatched element, of cardinality $mn$
(and thus dimension $mn-1$).
\end{itemize}
By Proposition~\ref{cor:f-is-Morse}, the matchings we obtain are Morse. Hence
Theorem~\ref{thm:open} follows from  Properties $(A)$ and $(B)$ using
Corollary~\ref{cor:Forman}. Let us now prove these properties. We
first study small values of $N$.
\begin{enumerate}
\item  If $N=1$, the graph is formed of isolated points. Hence the matching
$M(\emptyset, \emptyset, p)$ performed at the root of the tree is a
perfect matching. Property $(A)$  follows.
\item
For $N=2$, we leave it to the reader to check $(A)$ and $(B)$
for $M=1,2,3$, both for the graphs $\R(M,N)$ and $\sR(M,N)$. We then
proceed by induction on $M$, for $M\ge 4$.
\begin{enumerate}
  \item
For the
graph $\R(M,2)$, the root is a matching site, and its unique child is
$\Sig(A,B)$ with $A=\{(1,0)\}$ and $B=\{(0,0), (1,-1)\}$. The graph
obtained by deleting the vertices of $A\cup B$ is (a translate of)
$\sR(M-3,2)$. This shows that every unmatched element of
$\Sig(\R(M,2))$ is obtained by adding the vertex $(1,0)$ to (a
translate of) an unmatched element of $\Sig(\sR(M-3,2))$. The result
then follows by induction on $M$.
\item
For the graph $\sR(M,2)$, the root is  a splitting site. Its right
child is perfectly matched using the 
free pivot $(1,1)$, and has the empty
set as its unique child. The left child of the root is partially
matched using the pivot $(1,0)$. The unmatched elements are those of
$\Sig(A,B)$, with $A=\{(1,1)\}$ and $B=\{(1,0), (2,0)
\}$.  The graph
obtained by deleting the vertices of $A\cup B$ is (a translate of)
$\R(M-3,2)$. This shows that every unmatched element of
$\Sig(\sR(M,2))$ is obtained by adding the vertex $(1,1)$ to (a
translate of) an unmatched element of $\Sig(\R(M-3,2))$. The result
then follows by induction on $M$.
\end{enumerate}
\end{enumerate}

The study of the case $N=2$
leads us to introduce a notation that will
be useful in our forthcoming inductions.
\begin{Notation}
  Let $V$ be a subset of
vertices of the square grid, and assume that there exist
 $i,j\in \ZZ$ such that $V$ is the \emm  disjoint, union
$V=V_1\uplus \left( V_2+(i,j)\right)$ (where $V_2+(i,j)=\{v+(i,j):\
v\in V_2\}$). Let $X_1$ and $X_2$ be two collections of sets
 on the ground sets $V_1$ and $V_2$, respectively. Then
 $Y:=\{I_1\uplus \left( I_2+(i,j)\right): I_1\in X_1, I_2\in X_2\}$
is a collection of
 sets on the ground set $V$. We use the notation $Y\cong X_1 * X_2$
 to denote
 that the elements of $Y$ are formed by the concatenation of
 an  element of $X_1$ with (the translate of) an element of $X_2$. In
 particular, $|Y|= |X_1||X_2|$.
\end{Notation}
For instance, if $U(M,N)$ (resp.~$\sU(M,N)$) denotes the set of
unmatched elements in $\Sig(\R(M,N))$ (resp. $\Sig(\sR(M,N))$), the
above observations can be summarized by
$$
U(M,2)\cong U(3,2) * \sU(M-3,2) \quad \hbox{and}\quad
\sU(M,2)\cong \sU(3,2) * U(M-3,2).
$$
We  now return to our induction.
\begin{enumerate}\setcounter{enumi}{2}
\item The case $N=3$ is very similar to the case $N=2$. One first checks
 that the result holds for $M=1,2,3$, both for the graphs $\R(M,N)$
 and $\sR(M,N)$. For $M\ge 4$, the result is proved by induction on
 $M$, after observing that
$$
U(M,3)\cong U(3,3) * \sU(M-3,3) \quad \hbox{and}\quad
\sU(M,3)\cong \sU(3,3) * U(M-3,3).
$$
%
The following three observations will be useful in the rest of the proof.
Firstly, the tentative pivot is never taken in the third diagonal of slope
$-1$. 
Secondly, when it is taken in the second diagonal, it is a free
pivot.
Finally,  when the matching tree has a non-empty leaf $\Sig(A,B)$
(that is, when $M\not \equiv_3 1$), all vertices on the third diagonal
belong to $B$.
\item 
Now, let $N\ge 4$. The key  observation is that the top of the
matching tree coincides, as far as the prescribed vertices, pivots and
splitting sites are concerned, with the matching tree obtained for
$N=3$. This is illustrated in Figure~\ref{fig:match-open}, and holds
as long as the pivots are taken in the first two diagonals. Once these
pivots have been exhausted, we are left with (at most)  one non-empty
unmatched set $\Sigma(A,B)$, whose
prescribed vertices are those of $U(M,3)$ (if we work with $\R(M,N)$),
or  $\sU(M,3)$ (if we work with $\sR(M,N)$). Moreover, the tree
rooted at the  vertex $\Sigma(A,B)$ is isomorphic to the matching tree
of $\sR(M,N-3)$ (resp. $\R(M,N-3))$. This leads to 
\begin{align*}
 U(M,N)\cong U(M,3)* \sU(M,N-3),\\
\sU(M,N)\cong \sU(M,3)* U(M,N-3)
\end{align*}
when for $M\equiv_3 1$ both $U(M,3)$ and
$\sU(M,3)$ are empty. Properties $(A)$ and $(B)$ easily follow.
\end{enumerate}

\end{proof}

\section{Rectangles with cylindric boundary conditions}
\label{sec:cylindricrectangles}
We now study a ``cylindric'' version of the graphs $\R(M,N)$,
obtained by wrapping these graphs on a cylinder. For
$M, N \ge 0$  and $M$ \emm even,, we consider the graph
$\cR(M,N)$ obtained from $\R(M+1,N)$ by identifying the vertices
$(i,i)$ and $(M/2+i, -M/2+i)$, for $0\le i \le \lfloor
\frac{N-1}2\rfloor$.
Observe that the rectangles
$\sR(M,N)$ of Figure~\ref{fig:smoothrectangles}, when wrapped in a
natural way  around a cylinder, yield the same 
family  of graphs.
We denote by $\cZR(M,N)$ the alternating number of independent sets on 
the graph $\cR(M,N)$.

\begin{Theorem}\label{thm:cylindric}
 Let $M, N \ge1$, with $M$ even.  Denote $m= \lfloor
  \frac{M+1}3\rfloor$ and
$n= \lceil   N/3\rceil$.
 \begin{itemize}
    \item
  If $N\equiv_3 1 $, then $\Sigma(\cR(M,N))$ is
  contractible and $\cZR(M,N)=0$.

 \medskip
\item
Otherwise,
  \begin{itemize}
 \item If $M\equiv_3 0 $,  then $\Sigma(\cR(M,N))$ is homotopy equivalent
  to a wedge of $2^n$ spheres of dimension $mn-1$, and $\cZR(M,N)=2^n$.
    \item If $M\equiv_3 1 \mbox{ or }2 $ then $\Sigma(\cR(M,N))$
  is homotopy equivalent
  to a single sphere of dimension $mn-1$, and $\cZR(M,N)=(-1)^n$.
  \end{itemize}
 \end{itemize}
\end{Theorem}
\begin{proof}
We define a matching of $\Sig\equiv \Sigma(\cR(M,N))$ by adopting the
same choice of tentative pivots and splitting vertices as in the
proof of Theorem \ref{thm:open}.

\begin{figure}[thb]
\scalebox{0.9}{\input{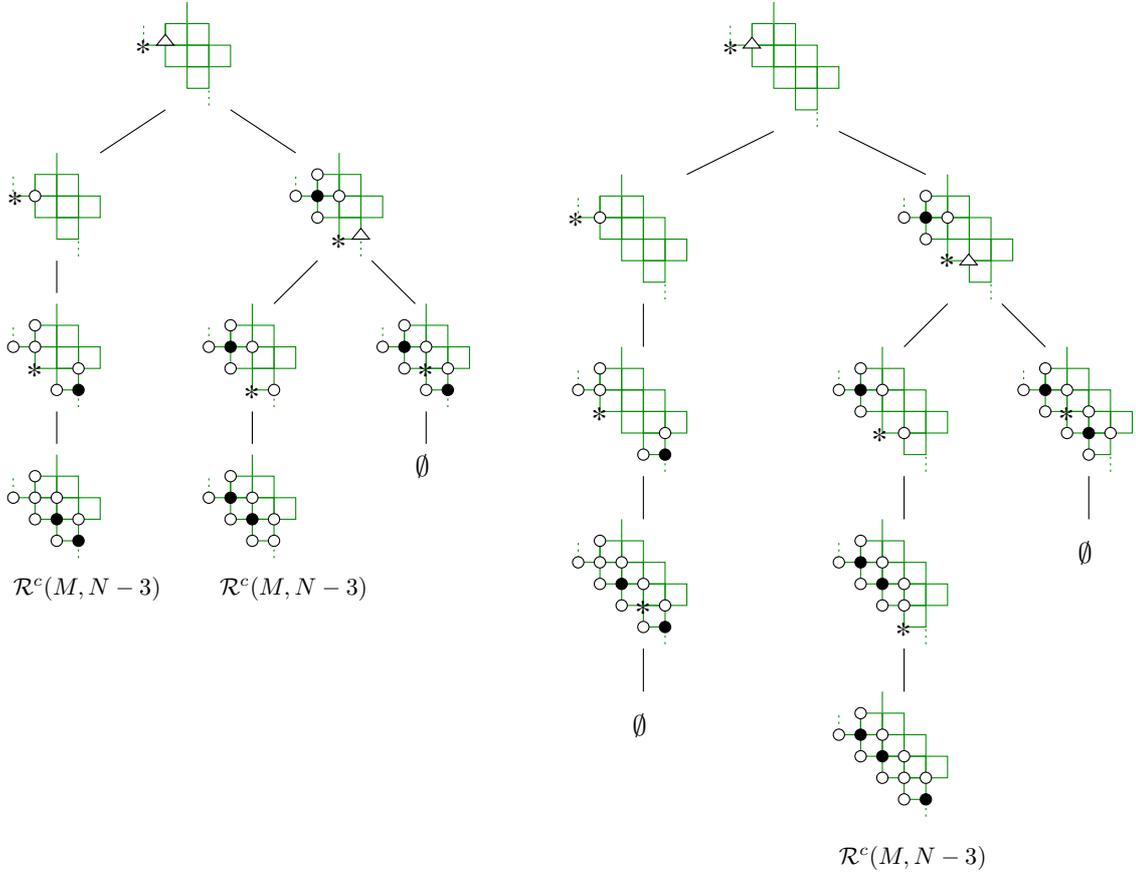}}
\caption{The top of the matching trees of $\Sigma(\cR(6,5))$
  and $\Sigma(\cR(8,5))$. This figure illustrates what happens for
  $M\equiv_3 0$ and $M\equiv_3 2$. We leave it to the reader to
  practice with the case $M\equiv_3 1$.}
\label{fig:match-cylindric}
\end{figure}

We are going to prove by induction on $N$ 
the following properties:
\begin{itemize}
  \item[$(A)$] if $N\equiv_3 1$, then there is no unmatched element
in $\Sig$,
\item[$(B)$] otherwise,
\begin{itemize}
 \item[$(B_1)$] If $M\equiv_3 0 $,  there are $2^n$  unmatched
  elements, each of    cardinality $mn$.
    \item[$(B_2)$] If $M\equiv_3 1 \mbox{ or }2 $ there is a unique unmatched
    element, of   cardinality $mn$.
  \end{itemize}
\end{itemize}
The theorem then follows from  Proposition~\ref{cor:f-is-Morse} and
Corollary~\ref{cor:Forman}.
Properties $(A)$ and $(B)$ are proved by induction on $N$.
\begin{enumerate}
\item When $N=1$, the graph is
formed  of isolated vertices, and the first pivot already matches
$\Sig(\cR(M,N))$ perfectly.
\item When $N=2$, the graph is a ring of $2M$ vertices.
Depending on the value of $M$ modulo 3, one finds two  (if $M\equiv_3
0$) or one  (if $M\equiv_3 1,2$) unmatched cell(s)
of cardinality $m$. These results are easily
obtained by constructing the whole matching tree of the graph.
We omit the details which are very similar to the proof of Theorem
\ref{thm:open}.
\item  The case $N=3$ is almost identical to $N=2$, since the pivot
  never appears 
in the third diagonal (Figure \ref{fig:match-cylindric}).
As in the case of rectangles
with open boundary conditions, for each
non-empty leaf $\Sig(A,B)$, all vertices on the rightmost diagonal
belong to $B$.  There are at most two such leaves (exactly two when
$M\equiv_3 0$). 
\item For $N\ge 4$,  the top of the tree
  coincides again with the matching tree of $\cR(M,3)$. Once the pivots of
  the first two diagonals have been exhausted, the prescribed
  vertices are exactly those of the
  first three diagonals. Denoting
 by $U^c(M,N)$ the set of unmatched elements of
  $\Sig(\cR(M,N))$, it follows that:
\begin{align*}
 U^c(M,N))\cong U^c(M,3)* U^c(M,N-3).
\end{align*}
\end{enumerate}
Properties $(A)$ and $(B)$ easily follow by induction on $N$.
\end{proof}

\section{Parallelograms with open boundary conditions}
\label{sec:paral}
\begin{figure}[thb]
\begin{center}
\input{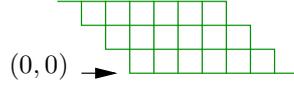}
\end{center}
\caption{The parallelogram graph $\Pa(4,8)$.}
\label{fig:paral}
\end{figure}

For $K, N\ge 1$, consider now the subgraph $\Pa(K,N)$ of the square
grid induced by the vertices $(x,y)$ satisfying
$$
0\le y \le K-1 \quad \mbox{ and }\quad -y\le x\le -y+N-1.
$$
An example is shown on Figure~\ref{fig:paral}. We denote by $\ZP(K,N)$
the alternating number of independent sets
on the graph $\Pa(K,N)$.

\begin{Theorem}\label{thm:open-paral}
 Let $K, N \ge1$. Denote $m= \lceil 2K/3\rceil$.
 \begin{itemize}
    \item
  If $K\equiv_3 1 $, then
  \begin{itemize}
    \item if $N\equiv_3 1 $ then $\Sigma(\Pa(K,N))$ is
  contractible and $\ZP(K,N)=0$,
\item otherwise,  $\Sigma(\Pa(K,N))$  is homotopy equivalent
  to a  sphere of dimension $mn-1$, with $n= \lceil N/3\rceil$,
and $\ZP(K,N)=(-1)^n$.
  \end{itemize}

 \medskip
\item  If $K\equiv_3 2 $,
write $N=2q K + r$, with $0\le r \le 2K-1$.
  \begin{itemize}
  \item
If $r\equiv_3 1,2$,  then $\Sigma(\Pa(K,N))$ is contractible and $\ZP(K,N)=0$,
\item otherwise $\Sigma(\Pa(K,N))$ is homotopy equivalent
  to a  sphere of dimension
$mn-1 $ with $n = \lceil  \frac {2K-1}{2K} \cdot \frac N 3\rceil$,
and $\ZP(K,N)=1.$
  \end{itemize}

 \medskip
\item  If $K\equiv_3 0 $,  write $N=2q(K+1)+r$ with $0 \le r  \le 2K+1$.
  \begin{itemize}
 \item If $r\equiv_3 0$ with $r\ge 1$, or $r\equiv_3 1 $ with
  $r\le 2K$,  then $\Sigma(\Pa(K,N))$ is
  contractible and $\ZP(K,N)=0$,
    \item otherwise, $\Sigma(\Pa(K,N))$  is homotopy equivalent
  to a  sphere of dimension $mn-1$
where  $n= \lceil \frac {2K+3}{2K+2} \cdot \frac N 3\rceil$,
and $\ZP(K,N)=1$.
  \end{itemize}
 \end{itemize}
\end{Theorem}
%
\noindent {\bf Remark.} 
We prove in Corollary~\ref{coro-0-paral} that for  $N\equiv_3 1$
  and $K$ large enough, the
  alternating number $\ZP(K,N;C,D)$ of independent sets on the
 parallelogram $\Pa(K,N)$ having prescribed conditions $C$ and $D$ on
  the top and bottom row is $0$, for all   configurations $C$ and
  $D$.
This is not in contradiction with the above theorem: for $K$ large
enough, if $K\not\equiv_3 1$,  the quotient $q$ appearing in the
theorem is simply $0$, so that the condition $N\equiv_3 1$ boils down
to $r\equiv_3 1$, with $r$ small compared to $K$, and $\ZP(K,N)=0$.
\begin{proof}
We construct a Morse matching of the graph $\Pa(K,N)$ by applying the
general method of Section~\ref{sec:preliminaries}. We then
 show that for all values of
$K$ and $N$, the matching thus obtained  has at most one unmatched
cell, of cardinality $mn$. As before, Corollary~\ref{cor:Forman}
completes the proof.

The results in the case $K\equiv _3 1$ are reminiscent of what we obtained for
rectangles (Theorem~\ref{thm:open}). 
The rule for choosing tentative pivots and
splitting vertices is the same as before, and the proof follows the same
principles. We do not repeat the argument.

The other two congruence classes of $K$ are more complicated, and 
require a  different pivot choice. Let us begin with the case $K\equiv_3
2$. First, we partition the set of vertices of $\Pa(K, \infty)$ into
triangular subsets $\T_1, \T_2, \ldots$ defined by:
 $$ \begin{array}{rcllrllllllllll} 
 \T_{2\ell +1}&= &\{(x,y) :& 0 \le
 y <K, & 2\ell K \le x+y , &x-y < 2\ell K+1\ \},\\ \T_{2\ell }&=&
 \{(x,y) :& 0 \le y <K, & 2(\ell -1)K+1 \le x-y, & x+y < 2\ell K\
\}.
\end{array}$$ 
These triangles are shown in Figure~\ref{fig:match-parallelogram} for
$K=5$.

\begin{figure}[tb]
\begin{center}
\scalebox{1.8}{\input{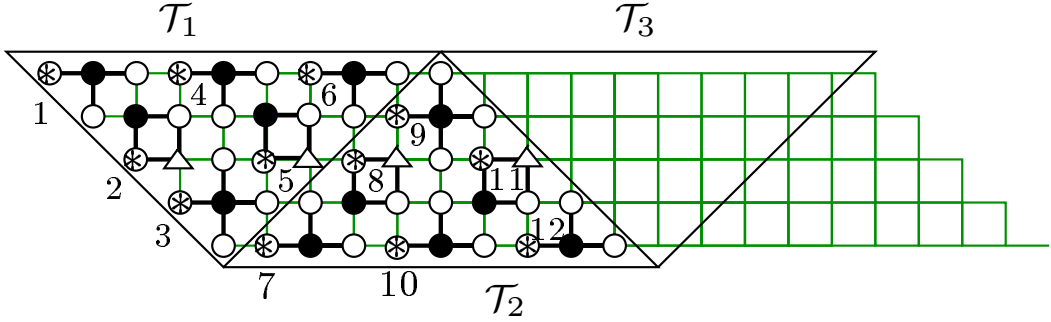}}
\end{center}
\caption{The configuration of prescribed vertices when the pivots of
  $\T_1$ and $\T_2$ have been exhausted, for $K=5$. This configuration
  occurs in the non-contractible branch of the matching tree. The pivots,
  indicated by $*$ and numbered by pivot order, and the splitting
  vertices $\triangle$, are
  empty. Each time a partial matching is performed, some new
  prescribed vertices appear: they are joined by thick line.}
\label{fig:match-parallelogram}
\end{figure}
We now describe the pivot order. 
The pivots 
are first  taken in $\T_{1}$, then in $\T_{2}$ (once all
vertices of $\T_{1}$ are prescribed),  and so on.  
The pivot order for $\T_{1}$ is very similar to what we have done previously.
That is, we follow the diagonals of slope $-1$ from upper left to lower right but with the
restriction that we must stay within the triangle. For $\T_{2}$ the
pivot order is what we get when turning $\T_{1}$ upside down. That is,
we follow the diagonals of slope 1 from bottom left to top right.

Denote by $U^p(K,N)$ the set of unmatched elements of
  $\Sig(\Pa(K,N))$. We observe (Figure~\ref{fig:match-parallelogram})
  that, once all pivots in  $\T_1$ and $\T_2$ have been exhausted,
  only one node $\Sig(A,B)$ of the matching tree is non-empty. The
  prescribed vertices of this node are those of  $A \cup B=\T_1 \cup
  \T_2$. This gives, for $N\ge 2K$:
$$
 U^p(K,N))\cong U^p(K,2K)* U^p(K,N-2K),
$$
which is the key for our induction on $N$. By convention, $\Pa(K,0)$
is the empty graph, and its unique independent set is the empty set,
so that $|U^p(K,0)|=1$. We also note that the
unique unmatched element of $\Sig(\Pa(K,2K))$ has cardinality
$m(m-1)=(2K+2)(2K-1)/9$.

Upon iterating the above identity, we obtain, if $N=2qK+r$,
$$
 U^p(K,N))\cong U^p(K,2qK)* U^p(K,r),
$$
where the only unmatched cell of  $\Sig(\Pa(K,2qK))$ has $qm(m-1)$
vertices. It thus remains to describe what our matching rule produces
for the graphs $\Pa(K,r)$, for $0\le r \le 2K-1$. 
The following properties are easily observed on the example of
Figure~\ref{fig:match-parallelogram}. 
\begin{enumerate}
  \item If $r\equiv _3 1$, the rightmost vertex 
  in the top row of $\Pa(K,r)$, which belongs to the triangle
  $\T_{1}$,  becomes a free pivot
  at some stage of the matching procedure,
  so  that $|U^p(K,r)|=0$.
 \item If $r\equiv _3 2$, 
the rightmost vertex in the top row of $\Pa(K,r)\cap \T_{2}$ becomes a free
pivot at some stage,
  so that again, $|U^p(K,r)|=0$.
\item Finally, if $r\equiv_3 0$, we obtain a unique unmatched cell, of
  cardinality $mr/3$.
\end{enumerate}
Putting together our recursion and the above results for $\Pa(K,r)$,
we find that the only non-contractible cases are when $r\equiv_3
0$. In this case, there is only one unmatched element in
$\Sig(\Pa(K,N))$, of cardinality $qm(m-1)+mr/3$. The result follows
for the case $K\equiv _3 2$.

\begin{figure}[thb]
\begin{center}
\scalebox{1.7}{\input{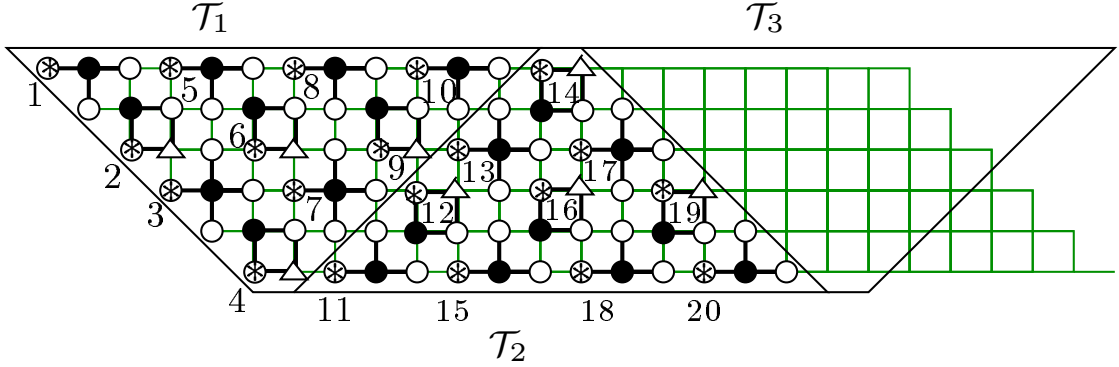}}
\end{center}
\caption{The configuration of prescribed vertices when the pivots of
  $\T_1$ and $\T_2$ have been exhausted, for $K=6$.}
\label{fig:match-parallelogram-bis}
\end{figure}

\medskip
Let us now adapt this to the final case $K\equiv _3
0$. The triangles that we used to define the pivot rule now become
trapezoids $\T_1, \T_2, \ldots$,  defined by: 
$$
\begin{array}{rcllrllllllllll}
\T_{2\ell +1}&= &\{(x,y) :& 0 \le y <K, & 2\ell(K+1) \le x+y ,
&x-y < 2\ell(K+1)+2\ \},\\
\T_{2\ell }&=& \{(x,y) :& 0 \le y <K, & 2(\ell-1)(K+1)+2 \le x-y,
& x+y < 2\ell( K+1)\  \}.
\end{array}
$$
These trapezoids are shown in Figure~\ref{fig:match-parallelogram-bis}
for $K=6$. 
The pivots are then chosen using the same rule as in the case
$K\equiv_3 2$.
Again, one observes that
when the pivots of  $\T_1$ and $ \T_2$ have been exhausted, only one
vertex $\Sig(A,B)$ of the matching tree is non-empty, and its
prescribed vertices are those of  $\T_1\cup \T_2$. 
This gives, for $N\ge 2K+2$:
$$
 U^p(K,N))\cong U^p(K,2K+2)* U^p(K,N-2K-2).
$$
Moreover, $\Sig(\Pa(K,2K+2))$ has a unique unmatched cell, of
cardinality $m(m+1)$. Let us write $N=2q(K+1)+r$, with $0\le r \le
2K+1$. Iterating the above identity gives
$$
 U^p(K,N))\cong U^p(K,2q(K+1))* U^p(K,r),
$$
where the only unmatched cell of $\Sig(\Pa(K,2q(K+1)))$ has
cardinality  $qm(m+1)$.

It  remains to describe what our matching rule produces
for the graphs $\Pa(K,r)$, for $0\le r \le 2K+1$. We refer again to
Figure~\ref{fig:match-parallelogram-bis}. 
\begin{enumerate}
  \item If $r\equiv _3 1$, with $r\le 2K$, 
  the rightmost vertex in the
  top row of   $\Pa(K,r)$ becomes a free pivot at some stage of the procedure, so
  that $|U^p(K,r)|=0$.
\item If $r=2K+1$, there is a 
unique unmatched cell, with cardinality
  $m(m+1)$ (it  coincides with the unmatched cell obtained for
  $\Pa(K,2K+2)$).
 \item If $r\equiv _3 0$ with $r>0$, 
the rightmost vertex in the top row of  $\Pa(K,r)\cap \T_2$ 
 becomes a free pivot at some stage of the matching procedure, so
  that again, $|U^p(K,r)|=0$.
\item If $r=0$, we have the empty graph, with the empty set as unique
  (and unmatched) independent set.
\item Finally, if $r\equiv_3 2$, we obtain a unique unmatched cell, of
  cardinality $m(r+1)/3$.
\end{enumerate}
Putting together our recursion and the above results for $\Pa(K,r)$,
we find that the only non-contractible cases are when $r=0$, $r=2K+1$
and $r\equiv_3 2$. In these cases, there is only one unmatched element in
$\Sig(\Pa(K,N))$, of cardinality $qm(m+1)+m\lceil r/3\rceil$. The
result follows for the case $K\equiv _3 0$.
\end{proof}

\noindent{\bf Remark.} The parallelogram  $\Pa(K,N)$ gives rise
to two distinct families of shapes with cylindric boundary conditions:
\begin{itemize}
  \item
Gluing the two diagonal borders of $\Pa(K,N+1)$ by identifying the
points $(-i,i)$ and $(-i+N,i)$ for $0\le i \le K-1$ gives the
``usual'' cylinder $\zs/N\zs \times \{0,1, \ldots, K-1\}$. 
It is conjectured in~\cite{JJ} that for odd $N$ the corresponding
alternating number of independent sets
 is 1, except for the case $N\equiv_6 3, K\equiv_3 1$ when it is
conjectured to be $-2$.

\item 
Gluing the two horizontal borders of $\Pa(K+1,N)$ by identifying the
points $(i,0)$ and $(i-K,K)$ for $0\le i \le N-1$ gives $\R^c(2K,N)$, the
cylindric version of the rectangle which we studied in
Section~\ref{sec:cylindricrectangles}.
\end{itemize}

\section{Transfer matrices}
\label{sec:transfer}

\subsection{Generalities}\label{sec:general}
We develop here a general (and very classical) transfer matrix
framework which we will instantiate later to the enumeration of
independent sets on various subgraphs of the square
lattice. See~\cite[Ch.~4]{stanley-vol1} for generalities on transfer
matrices.

Let $r\ge 1$, and let $\Ss$ be a collection of subsets in $\llb 1,
r\rrb:=\{1,2, \ldots, r\}$, of cardinality $d$.  Let $\TT$ be a square
matrix of size $d$, with complex coefficients,  whose rows and columns
are indexed by the elements of $\Ss$. The  entry of $\TT$
lying in row $C$ and column $D$ is denoted $\TT(C,D)$. Call
\emm configuration of length, $n$  any sequence $I=(C_0, C_1,
\ldots, C_n)$ of subsets of $\Ss$. We say that $C_0$ and
$C_n$ are the \emm borders, of $I$. The \emm weight, of $I$
is
$$
w(I)= \prod_{i=0}^{n-1} \TT(C_i,C_{i+1}).
$$
A \emm cyclic configuration of length, $n$ is a 
configuration $I=(C_0, C_1, \ldots, C_{n-1}, C_n)$
such that
$C_0=C_n$. Observe that for all $C \in \Ss$, $I=(C)$ is a cyclic
configuration of length $0$, so that there are exactly $d$ such configurations.

Let $t$ be an indeterminate.  It is well-known, and easy to prove,
that the length \gf\ of configurations with prescribed borders $C$ and
$D$, weighted as above, is
\beq\label{open-gf}
G_{C,D}(t):=\sum_{n\geq 0}\sum_{I=(C, C_1, \ldots, C_{n-1},D)} t^{n} w(I) =
(1-t\TT)^{-1}(C,D).
\eeq
From this, one derives that the length \gf\  of \emm cyclic,
configurations is the trace of $(1-t\TT)^{-1}$:
$$
G^{c}(t):=\sum_C G_{C,C}(t)= 
\rm{tr}(1-t\TT)^{-1}.
$$

In what follows, we will be interested in deriving eigenvalues of the
\emm transfer matrix, $\TT$  from the \gfs\ $G_{C,D}(t)$ and
$G^c(t)$. This is motivated by the conjectures of Fendley \emm et al.,
on the spectra of various transfer matrices related to the enumeration
of independent sets~\cite{Fendley}. One first observation is that
finding the whole
spectrum $(\lambda_1, \ldots, \lambda_{d})$ of the transfer matrix  is
equivalent to finding  the \gf\ of cyclic
configurations: indeed,
\beq\label{cylindric-gf}
G^{c}(t)= {\rm tr}(1-t\TT)^{-1}=\sum_{n\ge0} t^n {\rm tr}(\TT^n)
= \sum_{n\ge0} t^n(\la_1^n +\cdots +\lambda_{d}^n)
=\sum_{i=1}^d \frac 1{1-\lambda_it}.
\eeq
From~\eqref{open-gf},~\eqref{cylindric-gf}, 
and the classical formula  giving the inverse of a matrix in terms of
its determinant and its comatrix,
one concludes that
\begin{itemize}
  \item For every pair $C,D$, the reciprocals of the poles of $G_{C,D}(t)$ are
  eigenvalues of $\TT$, 
\item Conversely, the set of non-zero eigenvalues of $\TT$ coincides
  with the set   of reciprocals of poles
  of the series $G_{C,C}(t)$, for $C$ running over $\Ss$.
\end{itemize}
In other words: counting configurations with cyclic boundary
conditions gives the whole spectrum; at least partial information
can be derived from the enumeration of configurations with open
boundary conditions.

\begin{figure}[hbt]
\begin{center}
\input{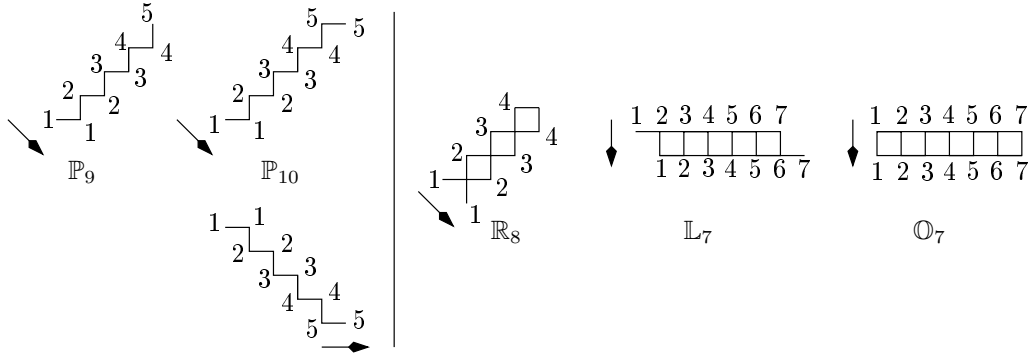}
\end{center}
\caption{How the various transfer matrices act.}
\label{fig:matrices}
\end{figure}

\smallskip
This general framework  will be specialized below to the case where the
sets $C$ and $D$ describe hard-particle configurations
(a.k.a. independent sets) on certain
layers on the square lattice.  The weights $\TT(C,D)$ will be designed
in such a way the weight of a configuration $I$ is $0$ if $I$ is not
an independent set, and  $(-1)^{|I|}$ otherwise.
 We have schematized in Figure~\ref{fig:matrices} the various transfer
 matrix we consider. They will be defined precisely in the text.

\subsection{Two complete spectra}
\label{sec:complete}

In this section, we combine the above generalities with the results
obtained in Section \ref{sec:cylindricrectangles} for the alternating number of
independent sets on the cylinders $\cR(M,N)$. As these cylinders can be
generated from two types of transfer matrices
(called $\RR_N$ and $\LL_N$ in Figure~\ref{fig:matrices}),
 we obtain the complete
spectra of these two families of transfer matrices. All the non-zero
eigenvalues are found to be roots of unity.

We begin with a transfer matrix  $\PP_N$ that describes the independent sets of
the 2-diagonal graph $\R(2,N)$. The rows of $\PP_N$ are
indexed by subsets $C$ of $\{1,2, \ldots, \lceil N/2\rceil\}$, its
columns are indexed by subsets $D$ of $\{1,2, \ldots, \lfloor
N/2\rfloor\}$, and
\beq\label{P-matrix-def}
\PP_N(C,D)=
\left\{
\begin{array}{ll}
i^{|C|+|D|}& \mbox{ if  }  C\cap D= C \cap (D+1)= \emptyset  ,\\
0 &\mbox{ otherwise.}
\end{array}\right.
\eeq
The notation $D+1$ means $\{i+1: i\in D\}$. 
Observe that $\PP_N$ is
not a square matrix if $N$ is odd\footnote{When $N$ is even, the
  matrix $\PP_N$ also describes the independent sets of the
  parallelogram $\Pa(N/2, 2)$, as illustrated in
  Figure~\ref{fig:matrices}. Accordingly, its powers will be used
  later to count independent sets on the parallelograms
  $\Pa(N/2,\cdot)$.}. 
However, if  $\bPP_N$ denotes the
transpose of $\PP_N$, then $\RR_N:=\PP_N\bPP_N$ is the (square)
transfer matrix corresponding to the
graph $\R (3,N)$: for $C, D \subseteq \{1, \ldots, \lceil
N/2\rceil\}$,
$$
\RR_N(C,D)=\sum_{E\subseteq \{1, \ldots, \lfloor N/2\rfloor\} }
\PP_N(C,E)\bPP_N(E,D)
= i^{|C|+|D|} \sum_{E : I=(C,E,D)\, {\rm     ind. set}} (-1)^{|E|}
$$
is, up to the factor $i^{|C|+|D|}$, the alternating number of
independent sets $I$ of $\R (N,3)$, with top and bottom borders $C$
and $D$ respectively.   The coefficients of this matrix have actually a
simpler expression. Indeed, the configuration $(C,E,D)$ is an
independent set if and only if $E$ is in the complement of $C\cup
D\cup (C-1) \cup (D-1)$. Hence the sum of terms $(-1)^{|E|}$ 
is 0 unless this complement is empty. That is, for $C, D
\subseteq\{1,2,\ldots , \lceil N/2\rceil\}$,
\beq\label{R-matrix-def}
\RR_N(C,D)=
\left\{
\begin{array}{ll}
i^{|C|+|D|}& \mbox{ if  }  \{1,2,\ldots , \lfloor N/2\rfloor\}=
C\cup D\cup (C-1) \cup (D-1)  ,\\
0 &\mbox{ otherwise.}
\end{array}\right.
\eeq
From~\eqref{P-matrix-def}, one derives that the entry $(C,C)$ in the product
$\RR_N^k=(\PP_N \bPP_N)^k$ is the 
alternating number of independent sets
on the cylinder $\cR(2k,N)$ studied in
Section~\ref{sec:cylindricrectangles},  with border condition $C$ on
the first diagonal. By Section~\ref{sec:general}, these numbers are
related  to the spectrum $(\lambda_1, \ldots, \lambda_d)$ of the
 transfer matrix $\RR_N$. More  precisely,  \eqref{cylindric-gf} gives:
\beq\label{trace-cylindric}
G^{c}(t)=d+\sum_{k\ge 1} \cZR(2k,N) t^k
=\sum_{i=1}^d \frac 1{1-\lambda_it},
\eeq
where the numbers $\cZR(2k,N)$ are given  in
 Theorem~\ref{thm:cylindric} and $d= 2^{\lceil  N/2\rceil}$.

\begin{Theorem}
Let $N\ge 1$. The  transfer matrix $\RR_N$ defined by {\rm
  \eqref{R-matrix-def}} has  size $2^{\lceil
  N/2\rceil}$. If $N\equiv_3 1$, then  $\RR_N$ is
  nilpotent (all its eigenvalues are 0).
Otherwise, $\RR_N$ has  eigenvalues:
\begin{itemize}
  \item $0$ with multiplicity $2^{\lceil
  N/2\rceil}-2^n$,
\item $1$ with multiplicity $(2^n+2(-1)^n)/3$,
\item $j$ and $j^2$ with multiplicity $(2^n-(-1)^n)/3$,
\end{itemize}
where $j=e^{2i\pi/3}$ and
$n=\lceil N/3\rceil$.
\end{Theorem}
\begin{proof}
  We start from the identity
$$
G^{c}(t)=d+\sum_{k\ge 1} \cZR(2k,N) t^k
={\rm tr} (1-t\RR_N)^{-1} = \sum_{i=1}^{d} \frac 1{1-\lambda_it}.
$$
This allows us to read off the
eigenvalues directly from the \gf\ of the numbers $ \cZR(2k,N)$, given
by Theorem~\ref{thm:cylindric}.

When $N\equiv_3 1 $, $G^c(t)=d$ and so all the eigenvalues of
$\RR_N$ are zero.
Otherwise,
$$
G^c(t)= d+ 2^n\,\frac{t^3}{1-t^3} + (-1)^n\frac{t+t^2}{1-t^3},
$$
and the result follows by a partial fraction expansion.
\end{proof}

In the case $N\equiv_3 1$, the above theorem gives an unexpected strengthening
of Theorem~\ref{thm:open}. This observation was communicated to us by Alan Sokal, merci \`a lui~! 
\begin{Corollary}\label{coro-0}
  Let $N\equiv_3 1$.  For $M> 2^{1+\lceil  N/2\rceil}$, the
  alternating number $\ZR(M,N;C,D)$ of independent sets on the
  rectangle $\R(M,N)$ having prescribed conditions $C$ and $D$ on
  the two extreme diagonals of slope $1$ is $0$, for all
  configurations $C$ and $D$.
\end{Corollary}
\begin{proof}
  First assume that $M=2m+1$, so that $m\ge 2^{\lceil  N/2\rceil}$. 
By the above theorem and the Cayley-Hamilton theorem, the $m$th
  power of $\RR_N$  vanishes. Thus by~\eqref{open-gf}, the series
  $G_{C,D} (t)$ is a polynomial in $t$ of degree at most $2^{\lceil
  N/2\rceil}-1$. But the coefficient of $t^m$ in this series is
  precisely $\ZR(2m+1,N;C,D)$.

Similarly, if $M=2m+2$ with $m\ge 2^{\lceil  N/2\rceil}$, the number
$\ZR(M,N;C,D)$ is the entry $(C,D)$ in the matrix $\RR_N^m \PP_N$,
which vanishes. 
\end{proof}

As observed at the end of Section~\ref{sec:paral},
 the cylindric shape
$\R^c(2k,N)$ can also be obtained by wrapping the parallelogram
$\Pa(k+1,N)$ on a cylinder, identifying the top and bottom
(horizontal) layers. Consequently, the results of Theorem~\ref{thm:cylindric} also give the
spectrum of another transfer matrix, denoted $\LL_N$, which describes
how to construct the shapes $\Pa(\cdot, N)$ layer by layer
(Figure~\ref{fig:matrices}). The rows and columns
of $\LL_N$ are indexed by independent sets of the segment
$\Pa(1,N)$. Thus the size of $\LL_N$ is the Fibonacci number $F_{N+1}$, with
$F_0=F_1=1$ and $F_{N+1}=F_N+F_{N-1}$, and if $C$ and $D$ are independent
sets of $\Pa(1,N)$,
\beq\label{LL-def}
\LL_N(C,D)=
\left\{
\begin{array}{ll}
i^{|C|+|D|}& \mbox{ if  }  C \cap (D+1)= \emptyset  ,\\
0 &\mbox{ otherwise.}
\end{array}\right.
\eeq
The generalities of Section~\ref{sec:general} imply that the spectrum
$(\mu_1, \ldots, \mu_d)$ of $\LL_N$ satisfies
$$
F_{N+1}+\sum_{k\ge 1} \cZR(2k,N) t^k
= \sum_{i=1}^{F_{N+1}} \frac 1{1-\mu_it}.
$$
Comparing with~\eqref{trace-cylindric} shows that  the spectra of
$\LL_N$ and $\RR_N$ coincide, apart from the multiplicity of the null
eigenvalue.
%
\begin{Theorem}
Let $N\ge 1$. The  transfer matrix $\LL_N$ defined by
{\rm\eqref{LL-def}} has  size
$F_{N+1}$. If $N\equiv_3 1$, then  $\LL_N$ is
nilpotent. Otherwise, $\LL_N$ has  eigenvalues:
\begin{itemize}
  \item $0$ with multiplicity $F_{N+1}-2^n$,
\item $1$ with multiplicity $(2^n+2(-1)^n)/3$,
\item $j$ and $j^2$ with multiplicity $(2^n-(-1)^n)/3$,
\end{itemize}
where $j=e^{2i\pi/3}$ and $n=\lceil N/3\rceil$.
\end{Theorem}
As for the matrix $\RR_N$, the nilpotent  case $N\equiv_3 1$ gives the
following corollary, which has to be compared  to
Theorem~\ref{thm:open-paral}.
\begin{Corollary}\label{coro-0-paral}
  Let $N\equiv_3 1$.  Then, for $K> F_{N+1}$, the
  alternating number $\ZP(K,N;C,D)$ of independent sets on the
 parallelogram $\Pa(K,N)$ having prescribed conditions $C$ and $D$ on
  the top and bottom row is $0$, for all
  configurations $C$ and $D$.
\end{Corollary}
\begin{proof}
  The number $\ZP(K,N;C,D)$ is the entry $(C,D)$ in the $(K-1)$th
  power of the transfer matrix $\LL_N$. But this power vanishes by the
  Cayley-Hamilton theorem.
\end{proof}
\subsection{Partial results on two other spectra}
\label{sec:partial}
We focus in this section on two transfer matrices that generate the
usual cylinder $\C(K,N):=  \{0,1,\ldots, K-1\}\times \zs/N\zs$.
This cylinder can be
obtained by identifying the diagonal borders  of the parallelogram
$\Pa(K,N+1)$. Alternatively, it can be obtained by wrapping the
ordinary $K\times (N+1)$ rectangle $\{0,1, \ldots, K-1\}\times \{0,1,
\ldots, N\}$ on a cylinder, identifying the vertices $(i,0)$ and
$(i,N)$.
Underlying the first construction are the matrices $\PP_{2K}$ defined
at the beginning of Section~\ref{sec:complete}
(see~\eqref{P-matrix-def} and Figure~\ref{fig:matrices}).
Underlying the second construction is  the  transfer matrix
$\OO_K$ that describes
how to construct the ordinary
rectangles of width $K$.
This matrix has size $F_{K+1}$, the $(K+1)$st Fibonacci number,
and its rows and columns are indexed by independent sets of the $K$
point segment. If $C$ and $D$ are two of these independent sets,
\beq
\label{O-matrix-def}
\OO_K(C,D)=
\left\{
\begin{array}{ll}
i^{|C|+|D|}& \mbox{ if  }  C \cap D= \emptyset  ,\\
0 &\mbox{ otherwise.}
\end{array}\right.
\eeq
The similarity with the definition~\eqref{LL-def} of the matrix
$\LL_N$ is striking, but the spectrum of $\OO_K$ is definitely more
complex than that of $\LL_N$. It is conjectured in~\cite{Fendley} that
all the eigenvalues of $\OO_K$ are roots of unity.

Let $\ZC(K,N)$ denote the alternating number of independent sets on the
cylinder $\C(K,N)$. From the generalities of
Section~\ref{sec:general}, we have:
$$
\sum_{N\ge 1} \ZC(K,N) t^N= {\rm tr}(1-t\PP_{2K})^{-1}-2^K =
{\rm tr}(1-t\OO_{K})^{-1}-F_{K+1}.
$$
That is, the spectra of the matrices $\PP_{2K}$ and $\OO_K$ coincide,
apart from the multiplicity of the null eigenvalue.

Alas, we do not know what the numbers $\ZC(K,N)$ are.
However, remember from Section~\ref{sec:general} that enumerative
results on configurations with open boundary conditions provide
partial informations on the spectrum of the transfer matrix.
Here, we 
 exploit the results of Section~\ref{sec:paral} on parallelograms
to obtain some information on the spectrum of $\PP_{2K}$ (and thus of
$\OO_{K}$).

For all $C,D\subseteq \{1, \ldots, K\}$, let $\ZP(K,N;C,D)$
be  the alternating number of independent sets of the parallelogram
$\Pa(K,N)$ having borders $C$ and $D$ respectively on the leftmost and
rightmost diagonal. Then~\eqref{open-gf} gives:
\beq\label{open-gf-1}
\begin{array}{lll}
  G_{C,D}(t)=
(1-t\PP_{2K})^{-1}(C,D)&=&
\displaystyle \delta_{C,D}
+(-i)^{|C|+|D|}
\sum_{I=(C, C_1, \ldots, C_{n-1},D), n\ge 1} (-1) ^{|I|}t^{n} \\
&=&\displaystyle\delta_{C,D}+
(-i)^{|C|+|D|}
\sum_{N\ge 1} \ZP(K,N+1;C,D)t^{N}.
\end{array}
\eeq
Note that in the first formula, the weight $(-i)^{|C|+|D|}(-1) ^{|I|}$
results in a weight  $i$ for each vertex of the extreme diagonals, as
it should.
Since the numbers $\ZP(K,N+2;\emptyset,\emptyset)$ coincide
with the numbers $\ZP(K,N)$ given in
Theorem~\ref{thm:open-paral}, this allows us to find some eigenvalues
of $\PP_{2K}$ and $\OO_{K}$. We indicate in the next section how our
pivot principle can be extended to count independent sets with
prescribed borders, so as to determine more series $G_{C,D}$ and thus
more eigenvalues of the matrix $\OO_K$.
\begin{Proposition}
\label{prop:partial-spectrum} The transfer matrix $\OO_K$ defined
by~{\rm\eqref{O-matrix-def}} satisfies the following properties.
\begin{itemize}
\item
If $K\equiv_3 1$, then $e^{i\pi/3}$ and $e^{-i\pi/3}$ are
eigenvalues of $\OO_K$.
\item If $K\equiv_3 2$, then all the $2K$th roots of unity, except
  maybe $-1$, are eigenvalues of  $\OO_K$.
\item If $K\equiv_3 0$, then all the $(2K+2)$th roots of unity, except
  maybe $-1$ and, if $K$ is odd, $\pm i$, are eigenvalues of  $\OO_K$.
\end{itemize}
\end{Proposition}
  \begin{proof}
Instantiating~\eqref{open-gf-1} to $C=D=\emptyset$ gives
$$
G_{\emptyset,\emptyset}(t)= (1-t\PP_{2K})^{-1}(\emptyset,\emptyset)
=1 +\sum_{N\ge 0}Z_\Pa(K,N) t^{N+1},
$$
with $\ZP(K,0)=1$. Recall that the reciprocals of the poles of this series are
eigenvalues of $\OO_K$ and $\PP_{2K}$. The numbers $\ZP(K,N)$ are
given in Theorem~\ref{thm:open-paral}. If $K\equiv_3 1$,
\beq\label{K1}
G_{\emptyset,\emptyset}=1+t+ \sum _{n \ge 1}(-1)^n \left( t^{3n}+
t^{3n+1}\right) = 
\frac 1 {1-t+t ^2}.
\eeq
If  $K\equiv_3 2$,
$$
G_{\emptyset,\emptyset}=1+ \sum _{q \ge 0}
\sum_{p=0}^{(2K-1)/3} t^{2qK+3p+1}
=1+ \frac t{1-t^{2K}} \frac{1-t^{2K+2}}{1-t^3}.
$$
Note that $({1-t^{2K+2}})/({1-t^3})$ is a polynomial, and that the
  only root this polynomial shares with  $1-t^{2K}$ is $-1$.

Finally, if  $K\equiv_3 0$,
\begin{multline*}
  G_{\emptyset,\emptyset}=1+ \sum _{q \ge 0}
t^{2q(K+1)+1} + \sum_{q\ge 0} t^{2q(K+1)+2K+2} +\sum_{q\ge 0}
\sum_{p=0}^{2K/3-1} t^{2q(K+1)+3p+3}\\
= \frac 1{1-t^{2K+2}}\left( 1+t + t^3\, \frac{1-t^{2K}}{1-t^3}\right).
\end{multline*}
Note that $( 1+t + t^3\,({1-t^{2K}})/({1-t^3}))$ is a polynomial. The
roots it shares with $(1-t^{2K+2})$ are $-1$, and, if $K$ is odd, $\pm
i$.
  \end{proof}

\section{Final comments and perspectives}
\subsection{Other quadrangles}\label{sec:quadrangles}
A natural generalization of the rectangles and parallelograms
studied in Sections~\ref{sec:rectangles} and~\ref{sec:paral}  is the
subgraph $\G(M,N)$ of the square lattice induced by the points $(x,y)$
satisfying
$$
ay\le x \le ay+M-1 \quad \mbox{ and } \quad  -bx\le y \le -bx+N-1,
$$
for given values of $a$ and $b$. We have solved above the cases $(a,b)=(1,1)$ and $(a,b)= (-1,0)$. In
particular, we have proved that in both cases, the alternating number of
independent sets is always $0$ or $\pm 1$. What about other values of
$a$ and $b$?
The case $(a,b)=(0,0)$,
which describes ordinary rectangles,
 shows that the simplicity of our results cannot
be extended to all pairs $(a,b)$. Indeed, even though the eigenvalues
of the transfer matrix $\OO_K$ are conjectured to be roots of unity,
the alternating  number $Z(K,N)$ of independent sets on a $K\times N$
rectangle does not show any obvious pattern. For instance, for $K=4$,
$$
\begin{array}{lll}
\displaystyle  \sum_{N\ge 0} Z(4,N) t^N &= &
\displaystyle\frac {1+{t}^{4}}{ \left( 1-t^2 \right)  \left( 1+t^3 \right) }\\
&=&1+{t}^{2}-{t}^{3}+2\,{t}^{4}-{t}^{5}+3\,{t}^{6}-2\,{t}^{7}+3\,{t}^{8}-3\,{t}^{9}+4\,{t}^{10}+O \left( {t}^{11} \right) .
\end{array}
$$
In particular, $Z(4,N)\sim(-1)^N N/3$ as $N $  goes to infinity.

In contrast, recall
 that it is  conjectured that for an ordinary rectangle with \emm cyclic,
 boundary conditions, the alternating number 
$\ZC(K,N)$
is~1 or  $-2$ when $N$ is odd~\cite{JJ}.

\begin{figure}[hbt]
\begin{center}
\input{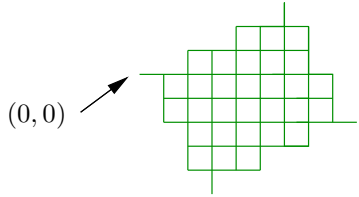}
\end{center}
\caption{The graph $\G(14,17)$ obtained for $a=b=2$.}
\label{fig:ab2}
\end{figure}

Still, the simplicity of the results obtained for rectangles and
parallelograms extends to other quadrangles. For instance,
if $a=b=2$, 
a pivot rule similar to the one used in the proof of
 Theorem~\ref{thm:open}
(the tentative pivot lies as high as possible on the leftmost line of
 slope $-2$)
 produces the following results, where the
alternating number of independent sets 
is denoted $Z_\G(M,N)$.
The proof is left to the reader.

\begin{Theorem}
  Let $M, N \ge 1$. Denote $m= \lceil M/5\rceil$ and  $n=  \lceil N/5\rceil$.
\begin{itemize}
  \item
If $M \equiv_5 0$ and $N\not = 3$, or if $M\equiv_5 1 $, or if $N \equiv_5  1,2$,
then $\Sigma(\G(M,N))$
is   contractible and $Z_\G(M,N)=0$.
\item
Otherwise,
   $\Sigma(\G(M,N))$ is homotopy equivalent
  to a sphere of dimension $mn-1$, and $Z_\G(M,N)=(-1)^{mn}$.
\end{itemize}
\end{Theorem}
It would be worth investigating which values of $a$ and $b$ produce
similar results.

\subsection{Transfer matrices}
 We have determined in Proposition~\ref{prop:partial-spectrum} some
 of the eigenvalues
of the ``hard'' transfer matrix $\OO_K$.  Recall that all its
 eigenvalues are conjectured to be roots of unity.
We give below the value of the characteristic polynomial $P(K)$ of the
matrix $\OO(K)$, for $1\le K\le 10$, and split this polynomial into
the part that is explained by Proposition~\ref{prop:partial-spectrum},
and the (bigger and bigger) part that is left unexplained. 
These data have been obtained with the help of Maple.
%
%
$$\begin{array}{rclcl}
  P(1) &=& \displaystyle \frac{1+t^3}{1+t} &\cdot& 1,\\
P(2)&=& \displaystyle \frac{1-t^4}{1+t} &\cdot& 1,\\
P(3)&=& \displaystyle \frac{1-t^8}{(1+t)(1+t^2)} &\cdot& 1,\\
P(4)&=&  \displaystyle \frac{1+t^3}{1+t} &\cdot& (1-t^2)(1-t ^4),\\
P(5)&=& \displaystyle \frac{1-t^{10}}{1+t} &\cdot& (1+t^4),\\
P(6)&=& \displaystyle \frac{1-t^{14}}{1+t} &\cdot& (1-t^4)^2,\\
P(7)&=& \displaystyle \frac{1+t^3}{1+t} &\cdot& \displaystyle \frac{(1+t^4)(1-t^{12})(1-t^{18})}{1+t^2},\\
P(8)&=& \displaystyle \frac{1-t^{16}}{1+t}&\cdot& (1-t^2)(1-t^4)^2(1+t^8)(1-t^{22}),\\
P(9)&=& \displaystyle \frac{1-t^{20}}{(1+t)(1+t^2)} &\cdot& \displaystyle \frac
{ (1+t^4)(1-t^{14}) (1+t^{10}) (1-t^{20}) (1-t^{26})}{1-t^2},\\
P(10)&=&\displaystyle \frac{1+t^3}{1+t} &\cdot& \displaystyle \frac
{(1-t^4)^2(1-t^{18})^2 (1-t^{24})^3 (1-t^{30})}{1+t^4}.
\end{array}$$
\medskip
By Section~\ref{sec:general}, we \emm know, that each missing factor
 occurs in at least one of the series $G_{C,C}(t)$ counting independent
 sets of the parallelogram $\Pa(K, \cdot)$ with
prescribed border $C$ on extreme diagonals. 
Conversely,
any  series  $G_{C,D}(t)$ may provide some of these missing factors
 (see~\eqref{open-gf}). Hence the following
question: can our pivot approach be recycled to compute
 some of these series, and  do we obtain new eigenvalues in this way?

For $C, D \subseteq \llb 1, K \rrb$, denote by
$\ZP (K,N;C,D)$ the alternating number
of independent sets of $\Pa(K,N)$ having border
conditions $C$ and $D$, respectively, on the first (last)
diagonal. (This notation was already introduced in
Section~\ref{sec:partial}.) Recall in 
particular the connection~\eqref{open-gf-1} between these numbers and the
series $G_{C,D}(t)$.

\begin{figure}[bht]
\begin{center}
\scalebox{1.7}{\input{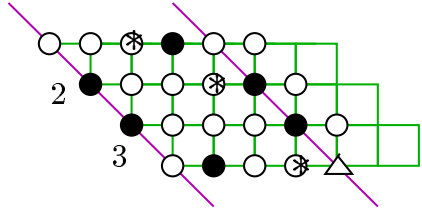}}
\end{center}
\caption{The pivot rule applied to $\Pa(4,N)$, with first diagonal
  $\{2,3\}$ and second diagonal $\{4\}$.}
\label{fig:paral-border}
\end{figure}

Take  $K=4$ and  $C=\{2,3\}$. If the configuration is $C$ on the first
diagonal, then, in the second diagonal, only the  vertex labeled 4 may
belong to an independent set. This gives
\beq\label{borders}
\begin{array}{lll}
\ZP (4,N; \{2,3\},D)&=&  \ZP (4,N-1; \emptyset,D) + \ZP(4,N-1; \{4\},D)\\
&=&\ZP (4,N-1; \emptyset,D) + \ZP(4,N-4; \{2,3\},D).
\end{array}
\eeq
The second identity is obtained by applying the pivot rule of
Section~\ref{sec:rectangles} to the independent sets counted by $\ZP (4,N-1;
\{4\},D)$ (Figure~\ref{fig:paral-border}). The above identity is 
valid for $N\ge 7$. We first 
specialize it to $D=\emptyset$, and work out what happens for small
values of $N$. Upon summing over $N$, we obtain
$$
G_{C,\emptyset}(t)= -\frac t{1-t^4}  \left( G_{\emptyset,
  \emptyset}(t)-t\right)= -\frac{t}{(1+t)(1-t+t^2)}.
$$
(We have used~\eqref{K1} for the value of $G_{\emptyset,
  \emptyset}$.)
We then specialize~\eqref{borders} to $D=C=\{2,3\}$. After working out
what happens for small values of $N$, we obtain
$$
G_{C,C}(t)=
\frac{1-tG_{\emptyset,C}(t)}{1-t^4}=\frac{1+t^2+t^3}{(1+t)(1-t^4)(1-t+t^2)},
$$
which now explains the missing factors $(1+t)(1-t^4)$ in $P(4)$.

For $K=5$, we obtain similarly the missing factor $1+t^4$ by
considering the numbers $\ZP(K,N; C,D)$, with $C=\{3,4\}$. Indeed, after a
discussion about the vertex labeled 5 in the second diagonal and a
few applications of the pivot rule, one finds
$$
\ZP (K,N;\{3,4\},D)=  \ZP (K,N-4; \emptyset,D) - \ZP(K,N-4;\{3,4\},D),
$$
from which we derive
$$G_{C,C}(t)=\frac 1{(1-t^5)(1+t^4)}.
$$
This gives the missing factor $(1+t^4)$.
\smallskip

It would be interesting to know  how far one can go with this
approach. That is, can we determine all series $G_{C,C}$ 
in this way?
This would allow us to count independent sets on the
ordinary cylinder, and thus to find all eigenvalues of the matrix
$\OO_K$.

Note that the largest cyclotomic factor that occurs in the polynomial
$P(K)$ seems to be $\Phi_{4K-10}$.

\bigskip
\noindent
\textbf{Remark added to the paper (07/03/2007)}: Sonja Cukic and
Alexander Engstr\"{o}m recently pointed out to
us that the following lemma~\cite[Lemma 2.4]{Engstrom} can be used to
derive some of our topological results.
\begin{Lemma}\label{lemma:Engstrom}
Let $v,w$ be vertices in a graph $G$. If $N(v) \subseteq N(w)$ then
$\Sigma(G)$ collapses onto $\Sigma(G-w)$. 
\end{Lemma}
Indeed, one can match all independent sets containing $w$ by adding or
removing the vertex $v$. 

For example, for the rectangles of Theorem \ref{thm:open}, one can
eliminate the vertices in every third diagonal  by a repeated
application of Lemma~\ref{lemma:Engstrom}. This leaves a graph
formed of several paths. One can then eliminate
every third vertex in each of these paths. The remaining graph is a
disjoint union of edges, and, if $M\equiv_3 1 $ or $N\equiv_3 1 $,
also isolated vertices. Thus $\Sigma(\R(M,N))$ collapses onto an
octahedral sphere (an $mn$-fold join of two points) or to the join of
an octahedral sphere with a simplex, respectively.
Theorem~\ref{thm:open} follows.

\bigskip \noindent
{\bf Acknowledgments:} We are grateful to Anders
       Bj\"{o}rner and Richard Stanley for inviting us to the
``Algebraic Combinatorics'' program at the Institut Mittag-Leffler
in Spring 2005, during which  part of this work was done.  All authors
were partially supported by the European Commission's IHRP
  Programme, grant HPRN-CT-2001-00272, ``Algebraic Combinatorics in
  Europe''.

\bibliographystyle{plain}
\bibliography{bibgrid}

\begin{thebibliography}{10}

\bibitem{Baxter}
R.~J. Baxter.
\newblock Hard hexagons: exact solution.
\newblock {\em J. Phys. A}, 13(3):L61--L70, 1980.

\bibitem{Engstrom}
A.~Engstr\"{o}m.
\newblock Independence complexes of claw-free graphs.
\newblock ArXiv:math.CO/0512420, 2005.

\bibitem{Fendley}
P.~Fendley, K.~Schoutens, and H.~van Eerten.
\newblock Hard squares with negative activity.
\newblock {\em J. Phys. A}, 38(2):315--322, 2005.
\newblock ArXiv:cond-mat/0408497.

\bibitem{Forman}
R.~Forman.
\newblock Morse theory for cell complexes.
\newblock {\em Adv. Math.}, 134(1):90--145, 1998.

\bibitem{FormanUser}
R.~Forman.
\newblock A user's guide to discrete {M}orse theory.
\newblock {\em S\'em. Lothar. Combin.}, 48:Art.\ B48c, 35 pp. (electronic),
  2002.

\bibitem{JJthesis}
J.~Jonsson.
\newblock {\em Simplicial complexes of graphs}.
\newblock PhD thesis, KTH, Stockholm, 2005.

\bibitem{JJdiag}
J.~Jonsson.
\newblock Hard squares on grids with diagonal boundary conditions.
\newblock {\em Preprint}, 2006.

\bibitem{JJ}
J.~Jonsson.
\newblock Hard squares with negative activity and rhombus tilings of the plane.
\newblock {\em Preprint}, 2006.

\bibitem{Munkres}
J.~R. Munkres.
\newblock {\em Elements of algebraic topology}.
\newblock Addison-Wesley Publishing Company, Menlo Park, CA, 1984.

\bibitem{sloane}
N.~J.~A. Sloane and S.~Plouffe.
\newblock {\em The encyclopedia of integer sequences}.
\newblock Academic Press Inc., San Diego, CA, 1995.
\newblock http://www.research.att.com/$\sim$njas/sequences/index.html.

\bibitem{stanley-vol1}
R.~P. Stanley.
\newblock {\em Enumerative combinatorics. {V}ol. $1$}, volume~49 of {\em
  Cambridge Studies in Advanced Mathematics}.
\newblock Cambridge University Press, Cambridge, 1997.

\end{thebibliography}

\end{document}